\documentclass[12pt]{article}
\usepackage{amsfonts}

\newcommand{\sect}[1]{\setcounter{equation}{0}\section{#1}}
\newcommand{\subsect}[1]{\subsection{#1}}

\renewcommand{\theequation}{\arabic{section}.\arabic{equation}}

\def\be{\begin{equation}}
\def\ee{\end{equation}}
\def\bea{\begin{eqnarray}}
\def\eea{\end{eqnarray}}

\def\1{\'{\i}}
\def\R{{\mathbb R}}

\def\Sch{{\cal S}}
\def\gal{\overline{{\cal G}}}

\def\a#1{a_{#1}}
\def\b#1{b_{#1}}
\def\c#1{c_{#1}}

\def\zt{h}
\def\zx{p}
\def\zv{k}
\def\zd{d}
\def\zc{c}
\def\zm{m}

\def\aa#1{\alpha_{#1}}

\def\ap{A_+}
\def\am{A_-}
\def\bp{B_+}
\def\bm{B_-}
\def\aaa{N}
\def\bb{M}

\def\masa{m}

\begin{document}

\thispagestyle{empty}

 \
\hfill\ J.Phys. A 33 (2000) 3445-3465

\
\vspace{1cm}

\begin{center}
 {\LARGE{\bf{(1+1)  Schr\"odinger Lie bialgebras}}}

 {\LARGE{\bf{and their Poisson--Lie groups}}}

\end{center}

\bigskip\bigskip

\begin{center} Angel Ballesteros, Francisco J. Herranz and Preeti
Parashar
\end{center}

\begin{center} {\it {Departamento de F\1sica, Universidad
de Burgos} \\   Pza. Misael Ba\~nuelos,
E-09001 Burgos, Spain}
\end{center}

\bigskip\bigskip

\begin{abstract}
All Lie bialgebra structures for the  $(1+1)$-dimensional
centrally   extended Schr\"o\-din\-ger algebra are explicitly
derived and proved to be of the coboundary type. Therefore,
since all of them come from a classical $r$-matrix, the complete
family of Schr\"odinger  Poisson--Lie groups can be deduced by
means of the Sklyanin bracket. All possible embeddings of the
harmonic oscillator, extended Galilei and $gl(2)$ Lie bialgebras
within the  Schr\"odinger classification are studied. As an
application, new quantum (Hopf algebra) deformations of the
Schr\"odinger algebra, including their corresponding quantum
universal
$R$-matrices, are constructed.
\end{abstract}

\newpage


\sect{Introduction}

The  $(1+1)$-dimensional centrally  extended  Schr\"odinger
algebra $\Sch$ is a non-semisimple Lie algebra spanned by  the
generators $\{D, H, P, K, C, M\}$ which obey  the following
commutation rules
\cite{Hagen,Nied}
\be
\begin{array}{lll}
[D, P]=-P \quad  &[D, K]=K\quad &[K, P]=M \cr
[D, H]=-2H \quad  &[D, C]=2C \quad &[H, C]= D
\cr
[K, H]=P \quad  &[K, C] =0 \quad
&[M,\,\cdot\,]=0\cr
[P, C]=-K\quad  &[P, H]=0\quad
& \cr
\end{array}
\label{aa}
\ee
where $D$ is the dilation, $H$ the  time translation, $P$ the
space translation, $K$ the Galilean boost, $C$ the conformal 
transformation, and $M$ is the mass (a central generator).
The Schr\"odinger algebra  contains many remarkable Lie
subalgebras: the harmonic oscillator algebra $h_4$ with
generators  $\{D, P, K, M\}$, the  $gl(2)$ algebra defined by
the  generators $\{D, H, C, M\}$, and the $(1+1)$ extended
Galilei algebra $\overline{{\cal G}}$ spanned by $\{H, P, K,
M\}$. In turn, the three dimensional Heisenberg--Weyl
algebra $h_3$ generated by $\{P, K, M\}$ is a subalgebra of
$h_4$; $gl(2)$ is isomorphic to a direct sum of $sl(2,\R)$, 
spanned by $\{D, H, C\}$, with the  central extension $M$; and
obviously,  the $(1+1)$  Galilei algebra ${\cal G}$ with
generators  $\{H, P, K\}$ is a subalgebra of $\gal$. Hence we
have the following embeddings:
\be
h_3 \subset  h_4 \subset \Sch  \qquad
sl(2,\R) \subset gl(2) \subset \Sch  \qquad
{\cal G} \subset \gal \subset \Sch .
\label{ab}
\ee

It is well-known that  $\Sch$  arises as the Lie-point  symmetry
algebra of  the  $(1+1)$ heat-Schr\"odinger equation (SE). This
relationship  has motivated the search for   $q$-deformed
analogues of $\Sch$  related to space-time discretizations  of
the  SE  on geometric lattices of the type
$x_n=q^n x_0$. In \cite{FV1,FV2} $q$-deformations of the vector
field realization of $\Sch$  were considered as symmetry
algebras of different discretized versions of the SE. From
another point of view, a different $q$-Schr\"odinger algebra was
directly introduced in \cite{Dobrev} in order to obtain a
generalized $q$-SE from its deformed representation theory.
However, none of these $q$-algebras has been found to be
consistent with a Hopf algebra structure. On the other hand, a
discretization of the SE on a regular space-time lattice
$x_n=x_0+n z$ has been also studied from a symmetry approach in
\cite{javier}. In this case, the resultant symmetry operators
close a non-deformed Schr\"odinger algebra.  Other $q$-SEs,
which are  related with quantum algebras different to  $\Sch$,
can be found in \cite{ita,DD,Micu}.

The first Hopf algebra deformations of $\Sch$ were introduced
in \cite{sch1,sch2}. These two quantum algebras were obtained
by starting from   (non-standard) quantum deformation of 
certain   subalgebras, namely, $h_4$
\cite{BH} and $gl(2)$ \cite{Aneva,BHN,Preeti}. Furthermore
these  quantum  $\Sch$ algebras lead to discretizations of the
SE on uniform lattices  \cite{sch3}. Recently it was shown in
\cite{BHNN} that, through a suitable non-linear change of basis,
the  algebra sector of the two  Schr\"odinger Hopf algebras
\cite{sch1,sch2}  can be mapped onto its classical counterpart,
thus preserving the deformation entirely in the coalgebra
sector. In this way  the discretized SEs of \cite{sch3} were
directly related  with the  ones previously obtained in
\cite{javier}.

 All these results  suggest the systematic  investigation of all
the possible quantum  Schr\"odinger algebras  that can be endowed
with a Hopf algebra structure. This problem leads,  as a
necessary first order approximation, to the obtention and
classification of all the Schr\"odinger Lie bialgebras.  This
is the main goal of the paper.   We obtain in Section 2 the
Lie bialgebra structures associated with the extended
Schr\"odinger algebra in $(1+1)$ dimensions. As $\Sch$ is a
non-semisimple Lie algebra our procedure will consist in
computing the most general cocommutator and, independently,
the most general classical $r$-matrix for $\Sch$; thus the
coboundary Schr\"odinger bialgebras can be identified. The final
result is that, despite the non-semisimple nature of $\Sch$, 
{\em all} its Lie bialgebra structures are coboundary ones, that
is,  they can be always obtained from a classical
$r$-matrix. This enables us  to construct their  corresponding
Poisson--Lie groups by making use of the Sklyanin bracket.

The following sections are devoted to show that  this
result contains a great amount of useful information
concerning quantum deformations of the Schr\"odinger
algebra. As a first step, in Section 3 we impose one 
generator to be primitive (besides $M$) in order to obtain
explicitly some relevant families of Lie bialgebras from the
previous classification. In particular, we are lead to three
different types of Lie bialgebras, that will be further
divided into two subfamilies each (standard and
non-standard).  Section 4 is devoted to the investigation of
the embeddings of  $h_4$, $gl(2)$ and $\gal$  Lie bialgebra
structures into the Schr\"odinger classification. This study
provides a systematic method to analyse the existence of the
extension of some known quantum deformations of all these
relevant subalgebras to the full Schr\"odinger algebra. 

As an application of this embedding technique, Section 5
presents two new examples of  $q$-deformed Schr\"odinger
algebras endowed with a Hopf algebra structure. In both
cases,  the universal quantum $R$-matrices can be also
obtained starting from the subalgebras. Finally, 
some remarks
concerning further applications of quantum Schr\"odinger
algebras are pointed out.


\sect{The  Schr\"odinger Lie bialgebras}

A {\em Lie bialgebra} 
$(g,\delta)$ is a Lie algebra $g$ endowed with a linear map
$\delta:g\to g\otimes g$  (the {\em cocommutator}) such that:

\noindent i) $\delta$ is a 1-cocycle, i.e.,
\be
\delta([X,Y])=[\delta(X),  1\otimes Y+ Y\otimes 1] +
[1\otimes X+ X\otimes 1,  \delta(Y)]   \qquad \forall X,Y\in g.
\label{ba}
\ee
\noindent ii) The dual map $\delta^\ast:g^\ast\otimes g^\ast \to
g^\ast$ is a Lie bracket on $g^\ast$.

 A Lie bialgebra $(g,\delta)$ is called a {\em coboundary} Lie
bialgebra if there exists an (skew-symmetric) element $r\in
g\otimes g$ (the {\em classical
$r$-matrix}) such that 
\be
 \delta(X)=[1\otimes X + X \otimes 1,\,  r]  \qquad  \forall
X\in g.
\label{bb}
\ee

Coboundary Lie bialgebras can be of two different types:

\noindent i)  {\em Non-standard} (or triangular): The $r$-matrix
is a skew-symmetric solution of the classical Yang--Baxter
equation (YBE):
\be
[[r,r]] = 0
\label{bc}
\ee
where $[[r,r]]$ is the Schouten bracket defined by
\be
[[r,r]] := [r_{12},r_{13}] + [r_{12},r_{23}] + [r_{13},r_{23}] .
\label{bd}
\ee
 If we denote $r = r^{ij} X_i \otimes X_j$, then $r_{12} =
r^{ij} X_i \otimes X_j \otimes 1, r_{13} = r^{ij}X_i \otimes 1
\otimes X_j$ and $r_{23} = r^{ij} 1 \otimes X_i \otimes X_j$.

 \noindent ii) {\em Standard} (or quasitriangular): The
$r$-matrix is a skew-symmetric solution of the modified
classical YBE:
\be
[X\otimes 1\otimes 1 + 1\otimes X\otimes 1 + 1\otimes 1\otimes X, \,
[[r,r]]\, ] = 0 \qquad \forall X \in g.
\label{be}
\ee

Finally, two Lie bialgebras $(g,\delta)$  and $(g,\delta')$ are
said to be {\em equivalent} if there exists an automorphism $O$
of $g$ such that
$\delta'=(O\otimes O)\circ\delta\circ O^{-1}$.


\subsect{The general solution}

The procedure  to characterize all Schr\"odinger Lie bialgebras
has two main steps. First it is necessary to deduce the most
general cocommutator, and second, to find out which of the
bialgebras so obtained come from a classical $r$-matrix. In
order not to burden the exposition, computational details  are
given in the appendix, and  the final result is summed up in the
following.

\noindent
{\bf Theorem 2.1.}  {\sl All $(1+1)$-dimensional centrally
extended Schr\"odinger Lie bialgebras are coboundary ones. They
are defined by  a classical  
$r$-matrix,  $r\in \Sch\wedge \Sch$, which   depends on 15 (real)
coefficients $a_i$, $b_i$, $c_j$ $(i=1,\dots,6;\ j=1,2,3)$:
\be 
\begin{array}{l}
r=\a1 D\wedge P + \a2 D\wedge H  +
\a3 P\wedge M + \a4 H\wedge M +
\a5 P\wedge H + \a6 P\wedge C\cr
 \quad\  + \b1 D\wedge K  + \b2 D\wedge C  +
\b3 K\wedge M    + \b4 C\wedge M +
\b5 K\wedge C + \b6 K\wedge H\cr
\quad\  +\c1 D\wedge M + \c2 P\wedge K +
\c3  H\wedge  C  
\end{array}
\label{ca}
\ee  
where the bialgebra coefficients  are subjected to 19 equations
casted into three sets:
\be 
\begin{array}{l}
 \a6^2 + \a6\b1 - 3\a1\b5 + \b5\b6=0\cr
 \a2\a3 - 2\a1\a4 - \a4\b6 - 3\a5\c1 + \a5\c2=0\cr
 \a1\a2 + \a2\b6 - \a5\c3=0\cr
 \a5\b1 - \a1\b6 - 2\a2\c1 - \a4\c3 =0\cr
  \a4\a6 + \a4\b1 - \a2\b3 - \a5\b4 + \a1\c1 - \a1\c2  =0\cr
3 \a1\b2 - \a2\b5 + \a6\c3 + \b1\c3 =0\cr
\a3\b2 + 2\a1\b4 - \a4\b5 + \a6\c1 + \a6\c2 + \b3\c3 =0\cr
3\a2\b5 + \b2\b6 - \a6\c3 =0
\end{array}
\label{cb}
\ee  
\be 
\begin{array}{l}
\b6^2 + \b6 \a1 - 3\b1\a5  +  \a5\a6  =0   \cr
\b2\b3- 2\b1\b4 - \b4\a6  - 3\b5\c1 - \b5\c2 =0\cr
\b1\b2 + \b2\a6 +  \b5\c3  =0\cr
 \b5\a1 -\b1\a6  - 2\b2\c1 + \b4\c3 =0\cr
\b4\b6 + \b4\a1    -\b2\a3  - \b5\a4 +  \b1\c1 + \b1\c2 =0\cr
3\b1\a2 - \b2\a5 - \b6\c3- \a1\c3   =0\cr
\b3 \a2 + 2\b1\a4  - \b4\a5  + \b6\c1 - \b6\c2 - \a3\c3  =0\cr
3\b2\a5 + \a2\a6 +  \b6\c3 =0 
\end{array}
\label{cc}
\ee  
\be 
\begin{array}{l}
4\a2\b2 + \c3^2 =0\cr
2\a4\b2 + 2\a2\b4 + \a5\b5 - \a6\b6 =0\cr
 2\a1\b1 -\a1\a6  - \b1\b6 + \a5\b5 - \a6\b6 =0 .
\end{array}
\label{cd}
\ee
The Schouten bracket of $r$ (\ref{ca}) turns out to be
\be
[[r,r]]= (\a3\a6+ \b3\b6 - \a3\b1 - \a1\b3   - \c2^2)K\wedge
M\wedge P 
\label{ce}
\ee
and the two types of classical $r$-matrices are distinguished
by an additional equation:
\be 
\begin{array}{ll}
\mbox{Standard:}&\quad \a3\a6+ \b3\b6 -  \a3\b1 - \a1\b3  -
\c2^2\ne 0\cr
 \mbox{Non-standard:}&\quad \a3\a6+ \b3\b6 - \a3\b1 - \a1\b3  -
\c2^2= 0 .
\end{array}
\label{cf}
\ee}

\medskip

The classical $r$-matrix (\ref{ca}) gives rise to the  general
cocommutators by applying the relation (\ref{bb}):
\be 
\begin{array}{l}
 \delta(D)= -\a1 D\wedge P -2\a2 D\wedge H-\a3 P\wedge M 
 - 2 \a4 H\wedge
M - 3 \a5 P\wedge H + \a6 P\wedge C\cr
\qquad\quad 
+\b1 D\wedge K +2\b2 D\wedge C+\b3 K\wedge M + 2 \b4 C\wedge M
+ 3 \b5 K\wedge C- \b6 K\wedge H\\[6pt]
 \delta(P)=(\a6-\b1)K\wedge P - \a2 H \wedge P  - \b1 D\wedge M
-\b2(C\wedge P + D\wedge K) \cr
\qquad\quad  -\b4 K\wedge M + \b5 C\wedge M + \b6 H\wedge
M + (\c1 - \c2) P\wedge M + \c3 K\wedge H\\[6pt]
\delta(K)=(\a1-\b6) P\wedge K + \b2 C \wedge K +  \a1 D\wedge M
+\a2(H\wedge K + D\wedge P)\cr
\qquad\quad   +\a4 P\wedge M - \a5 H\wedge M - \a6 C\wedge
M  - (\c1 + \c2) K\wedge M + \c3 P\wedge C\\[6pt]
\delta(H)=
 - (2\a1 + \b6) P\wedge H 
- (\a6 + \b1) D\wedge P
- 2 \b1 K\wedge H - 2 \b2
C\wedge H\cr
\qquad\quad - \b3 P\wedge M + \b4 D\wedge M    -  \b5 (D\wedge K
+ P\wedge C) + 2 \c1 H\wedge M -\c3 D\wedge H\\[6pt]
\delta(C)=  (2\b1 + \a6) K\wedge C+   (\b6 + \a1) D\wedge K + 2
\a1 P\wedge C + 2 \a2 H\wedge C\cr
\qquad\quad + \a3 K\wedge M - \a4 D\wedge M  + \a5 (D\wedge P +
K\wedge H) - 2 \c1 C\wedge M  -\c3 D\wedge C\\[6pt]
\delta(M)=0.
\end{array}
\label{ci}
\ee

 The most general element $\eta\in \Sch \otimes \Sch$ which is
$Ad^{\otimes 2}$ invariant is simply $\eta=\tau M\otimes M$
where $\tau$ is an arbitrary real number. The classical
$r$-matrix $r'=r+\eta$ gives rise to the same Schr\"odinger
Lie bialgebra than $r$, that is, (\ref{ci}); hence $r'$ is
the most general  non-skewsymmetric classical $r$-matrix for
$\Sch$.

On the other hand, the Schr\"odinger algebra automorphism 
  defined by
\be 
\begin{array}{lll}
D\to -D&\qquad P\to -K &\qquad K\to -P\cr
M\to -M&\qquad H\to -C&\qquad C\to -H  
\end{array}
\label{cg}
\ee
(which leaves the Lie brackets (\ref{aa}) invariant)  can be
 implemented at a bialgebra level by introducing a suitable
transformation of the  parameters $a_i$, $b_i$ and $c_j$ given by
\be 
\begin{array}{lll}
a_i\to b_i&\qquad b_i\to a_i&\qquad i=1,\dots,6\cr
\c1\to \c1&\qquad \c2\to -\c2&\qquad\c3\to -\c3 .
\end{array}
\label{ch}
\ee
We stress that under the bialgebra automorphism defined by the 
maps (\ref{cg}) and (\ref{ch}),   the general classical
$r$-matrix (\ref{ca}), the equations  (\ref{cd}) and the
Schouten bracket (\ref{ce}) remain invariant, while  the sets of
equations (\ref{cb}) and (\ref{cc}) are interchanged. As
expected, we also find that $\delta(P) \leftrightarrow
\delta(K)$, $\delta(H) \leftrightarrow \delta(C)$,  while
 $\delta(D)$ and $\delta(M)$ remain unchanged.


\subsect{The  Schr\"odinger Poisson--Lie groups}

Since all Schr\"odinger Lie bialgebras  are coboundary ones, we
can obtain their corresponding Poisson--Lie groups by means of
the Sklyanin bracket provided by a given classical $r$-matrix 
$r=\sum_{i,j} r^{ij}X_i\otimes X_j$ \cite{Drlb}: 
\be
\{f,g\}=\sum_{i,j} r^{ij}(X_i^Lf\, X_j^L g -
X_i^Rf\, X_j^R g)
\label{ma}
\ee
where   $X_i^L$ and $X_j^R$  are  left and right invariant vector
fields on the Schr\"odinger group. Thus   we consider the
following $4\times 4$ real matrix representation of  the
centrally extended Schr\"odinger algebra \cite{Patera}:
\be
\begin{array}{l}
H=\left(\begin{array}{cccc}
 0 &0 & 0 &0  \\ 0 & 0 &  1&0 \\ 0 & 0 & 0 & 0\\ 0 & 0 & 0 & 0
\end{array}\right)\quad
 P =\left(\begin{array}{cccc}
 0 & 0 & 1 &0  \\ 0 & 0 &  0&1 \\ 0 & 0 & 0 &0\\ 0 & 0 & 0 & 0
\end{array}\right)\quad 
K =\left(\begin{array}{cccc}
 0 & 1  &0 &0 \\ 0 & 0 &  0&
0   \\0 & 0 & 0 & -1\\0 & 0 & 0 & 0
\end{array}\right)\cr
\cr
D =\left(\begin{array}{cccc}
 0 &0 &0 &0 \\ 0 & -1 &  0&
0   \\ 0 &0 & 1 & 0\\ 0 & 0 & 0 & 0
\end{array}\right)\quad
 C =\left(\begin{array}{cccc}
 0 & 0 & 0 &0  \\ 0 & 0 &  0 & 0 \\ 0 &-1 & 0 & 0\\ 0 & 0 & 0 &
0 \end{array}\right)\quad
M =\left(\begin{array}{cccc}
 0 & 0 & 0 &2  \\ 0 & 0 &  0&
0   \\ 0 & 0 & 0 & 0\\ 0 & 0 & 0 & 0
\end{array}\right)  .
\end{array}
\label{mb}
\ee
This allows us to write an element of the  extended
Schr\"odinger group as
\bea
 &&g=\exp\{\zm M\}\exp\{\zx P\}\exp\{\zv K\}\exp\{\zt
H\}\exp\{\zc C\}\exp\{\zd D\}\nonumber \\[4pt]
&&\quad=\left(\begin{array}{cccc}
  1 & \bigl(\zv-(\zx+\zv\zt)\zc\bigr)  e^{-\zd}& (\zx+\zv
\zt)e^\zd &2\zm-\zx
\zv  \\  
0 &(1-\zt\zc)e^{-\zd} & \zt e^\zd& \zx \\ 
0 &-\zc e^{-\zd} & e^\zd &-\zv\\ 0 & 0 &0& 
1\end{array}\right) .\label{mc}
\eea
Then left and right invariant vector
fields can be  deduced
\bea
&&X^L_H=  e^{-2\zd} \left\{ {\partial_\zt}  - \zc 
{\partial_\zd}-\zc^2
 {\partial_\zc}\right\}  \cr
 &&X^L_P=(1-\zc\zt) e^{-\zd} \left\{ {\partial_\zx} +  \zv  
\partial_\zm
\right\} + \zc
e^{-\zd}{\partial_\zv} \cr
&&X^L_K=e^{\zd}\left\{ {\partial_\zv} -\zt   {\partial_\zx}   -
\zt \zv   \partial_\zm\right\} \cr
&&X^L_D=\partial_\zd \qquad X^L_C=e^{2\zd}\partial_\zc \qquad
X^L_M=\partial_\zm 
\label{md}
\eea
\bea
 &&X^R_H= {\partial_\zt}-\zv{\partial_\zx} - \frac{1}2
\zv^2\partial_\zm 
\qquad X^R_P=\partial_\zx\cr
&&X^R_K= {\partial_\zv}+\zx \partial_\zm 
\qquad X^R_M=\partial_\zm\cr
&&X^R_D= {\partial_\zd}+2\zc{\partial_\zc}-2\zt\partial_\zt
-\zx\partial_\zx + \zv \partial_\zv \cr
 && X^R_C=-\zt\partial_\zd
+(1-2\zt\zc)\partial_\zc+\zt^2\partial_\zt +\zx\partial_\zv+
\frac{1}2 \zx^2\partial_\zm .
\label{me} 
\eea
 We substitute now these vector fields  and the $r$-matrix
(\ref{ca}) within the  Sklyanin bracket (\ref{ma}) and compute
the Poisson--Lie brackets amongst the group coordinates
$\{\zd,\zt,\zx,\zv,\zc,\zm\}$; they turn out to be
\bea
 &&\{\zd,\zt\}=  \a2 ( e^{-2\zd}-1) + \b2 \zt^2 -\c3 \zt
\nonumber\\[2pt] &&\{\zd,\zx\}=\a1 \left( e^{-\zd} (1-\zc\zt)
-1\right)+\a2\zv+\a5 \zc e^{-3\zd} (1-\zc\zt) -\a6 \zt\cr
&&\qquad\qquad\quad -\b1 \zt e^{\zd} + \b2 \zt\zx - \b6 \zt \zc
e^{-\zd}+\c3 \zt\zv \nonumber\\[2pt]
 &&\{\zd,\zv\}=\a1 \zc e^{-\zd} + \a5 \zc^2 e^{-3\zd} +\b1 (
e^{\zd}-1) -\b2 (\zx +\zt\zv) -\b5 \zt +\b6 \zc e^{-\zd}
\nonumber\\[2pt] &&\{\zd,\zc\}=-\a2 \zc^2 e^{-2\zd}+\b2
(e^{2\zd}-1) -\c3\zc\nonumber\\[2pt] 
&&\{\zt,\zx\}= 2 \a1 \zt - \a2 (\zx + 2 \zt\zv) +
\a5 \left( 1 -  e^{-3\zd}(1-\zc\zt)\right)+ \a6 \zt^2\cr
&&\qquad\qquad\quad  - \b2 \zx\zt^2+\b6 \zt e^{-\zd} -\c3
\zv \zt^2 \nonumber\\[2pt]
 &&\{\zt,\zv\}=\a2 \zv- \a5 \zc e^{-3\zd} + 2 \b1 \zt +\b2 \zt
(2\zx  +\zt\zv)+\b5 \zt^2 + \b6 ( 1-e^{-\zd}) -\c3 \zx
\nonumber\\[2pt] &&\{\zt,\zc\}=2\a2 \zc +2\b2 \zt(1-\zc\zt) +
2\c3 \zt \zc \nonumber\\[2pt] &&\{\zx,\zv\}=\a1 \zv -\a2\zv^2
-\a6 \zx+\b1 \zx+\b2 \zx^2-\b6 \zv +\c3
\zv\zx \nonumber\\[2pt]
 &&\{\zx,\zc\}=2\a1\zc - 2 \a2 \zv\zc- \a5 \zc^2 e^{-3\zd}
(1-\zc\zt)+\a6
\left(e^{\zd}(1-\zc\zt)-1 + 2\zc\zt \right)\cr
&&\qquad\qquad\quad +\b2 \zx (1-2\zc\zt)-\b5 \zt e^{3\zd} +\b6
\zt\zc^2 e^{-\zd} +\c3 \zv(1-2\zc\zt) \nonumber\\[2pt]
&&\{\zv,\zc\}=-\a5  \zc^3 e^{-3\zd} +\a6\zc  e^{\zd} +2\b1 \zc
-\b2 \left( \zv (1-2\zc\zt) - 2\zc\zx\right)\cr
&&\qquad\qquad\quad 
+\b5 (e^{3\zd}-1 + 2 \zc\zt)-\b6   \zc^2 e^{-\zd}\label{mf}\\
 &&\{\zm,\zd\}=-\a1  \zv  e^{-\zd} (1-\zc\zt)-\frac 12 \a2
\zv^2+\a4 
\zc e^{-2\zd} -\a5 \zv  \zc e^{-3\zd} (1-\zc\zt)\cr
 &&\qquad\qquad\quad +\b1 (\zx +  \zt\zv e^{\zd})+\frac 12 \b2
\zx^2-\b4
\zt+\b5 \zt \zx +\b6 \zt\zv \zc e^{-\zd}  -\frac 12 \c3 \zt
\zv^2\nonumber\\[2pt]
 &&\{\zm,\zt\}= \a2 \zt\zv^2+\a4 (1-e^{-2\zd})+\a5
\zv e^{-3\zd} (1-\zc\zt)-2\b1\zt\zx-\b2 \zt\zx^2\cr
 &&\qquad\qquad\quad +\b4 \zt^2-\b5 \zx\zt^2-\b6 (\zx + \zt\zv
e^{-\zd})-2
\c1\zt +\frac 12 \c3 (\zx^2 +\zt^2\zv^2) \nonumber\\[2pt]
 &&\{\zm,\zx\}= \frac 12 \a2 \zx\zv^2 +\a3
\left(1-e^{-\zd}(1-\zc\zt)
\right) -\a4 \zv-\frac 12 \a5 \zv^2+\frac 12 \a6 \zx^2\cr
&&\qquad\qquad\quad -\b1 \zx^2 -\frac 12 \b2
\zx^3+\b3 \zt e^{\zd} +\b6 \zv\zx-\c1 \zx +\c2 \zx -\frac 12\c3
\zv\zx^2\nonumber\\[2pt]
 &&\{\zm,\zv\}=-\frac 12 \a2 \zv^3-\a3 \zc e^{-\zd}  +\b1
\zv\zx+\frac 12
\b2
\zv\zx^2-\b3 (e^{\zd}-1)+\b4 \zx\cr
&&\qquad\qquad\quad
-\frac 12 \b5 \zx^2-\frac 12 \b6 \zv^2
+\c1 \zv +\c2 \zv +\frac 12 \c3 \zx\zv^2\nonumber\\[2pt]
 &&\{\zm,\zc\}=-\a2 \zc\zv^2+\a4 \zc^2 e^{-2\zd} -\a5 \zv \zc^2
e^{-3\zd} (1-\zc\zt)+\a6 \zv e^{\zd} (1-\zc\zt)\cr
&&\qquad\qquad\quad
+2\b1\zc\zx+\b2 \zc\zx^2-\b4
 (e^{2\zd}-1 + 2 \zc\zt)-\b5 \left( \zt \zv e^{3\zd} + \zx
(1-2\zc\zt)
\right) \cr
&&\qquad\qquad\quad
 +\b6  \zt\zv \zc^2  e^{-\zd} + 2\c1 \zc +\frac 12 \c3 \zv^2 (1
- 2 \zc\zt) .
\nonumber
\eea 
 The parameters $a_i$, $b_i$ and $c_j$ must satisfy the 19
relations (\ref{cb})--(\ref{cd}); it is a matter of cumbersome
computations to recover these conditions  from the Jacobi
identities  coming from the Poisson algebra (\ref{mf}). On the
other hand, it can be checked that, as expected, the linear
terms of these Poisson--Lie brackets lead  to   the dual of the
cocommutators  (\ref{ci}).


\sect{The Schr\"odinger Lie bialgebras with two\\ primitive
generators}

A common feature of all extended Schr\"odinger bialgebras is
that the central generator $M$ is always  {\em primitive},
i.e., its cocommutator vanishes $\delta(M)=0$; this, in turn,
means that   its coproduct  can be always taken as
$\Delta(M)=1\otimes M+M\otimes 1$ after deformation. In
general, the existence of  primitive generators strongly
determines as much mathematical as possible   physical
properties of the corresponding quantum deformation  (see for
example
\cite{sch1,sch2}). Therefore, in order to specialize our
general results into several multiparametric families and to
obtain explicit quantization results, we  impose {\em one}
additional generator
$X$ to be primitive (besides $M$). Furthermore, the condition
$\delta(X)= 0$ rather simplifies the   equations
(\ref{cb})--(\ref{cf}).  Finally, due to the automorphisms
(\ref{cg}) and (\ref{ch}), there is no loss of generality if
we restrict our study to three types of bialgebras: those
with either $D$, $P$ or $H$ primitive.


\subsect{D primitive}

 The condition $\delta(D) = 0$  leaves $\c1$, $\c2$ and $\c3$ 
as the initial free  parameters, all the remaining ones being
equal to zero.  The equations (\ref{cb})--(\ref{cd}) imply that
$\c3 = 0$, so that   the Schouten bracket reduces  to $[[r,
r]] = - \c2^2 K \wedge M \wedge P$. Hence, the condition
$\c2\ne 0$ leads to   the standard subfamily of Lie
bialgebras, while $\c2=0$ gives the non-standard one.
Therefore each subfamily of Lie bialgebras is characterized
by the following parameters:

\noindent
$\bullet$ Standard subfamily:  $\c1$, $\c2 \ne
0$.

\noindent
$\bullet$ Non-standard subfamily: $\c1$.

From now on, it will be understood that all the parameters 
that do not appear explicitly for each subfamily are equal to
zero. We will also write  the  constraints that the
non-vanishing parameters must satisfy.

 The resulting $r$-matrix, cocommutators and non-vanishing
Poisson--Lie brackets for these subfamilies are given by
\bea 
&&r=\c1 D\wedge M+\c2   P\wedge K  \nonumber\\[2pt]
&&\delta(D)=0\qquad \delta(M)=0\label{na} \\
&&\delta(P)=(\c1-\c2) P \wedge M\qquad 
\delta(K)=-(\c1+\c2) K\wedge M\nonumber\\[2pt]
&&\delta(H)=2\c1 H\wedge M\qquad 
\delta(C)=-2\c1  C\wedge M .
\nonumber
\eea
\be 
 \{\zm,\zt\}= -2
\c1\zt \qquad
 \{\zm,\zx\}=  -(\c1-\c2 ) \zx  
\qquad
 \{\zm,\zv\}= (\c1+\c2) \zv \qquad
\{\zm,\zc\}= 2\c1 \zc .
\label{nb}
\ee  
 Obviously, the non-standard Lie bialgebra corresponds to the
substitution
$\c2=0$ in these expressions.


\subsect{P primitive}

 If we impose now $\delta(P) = 0$, we find that the   initial
parameters  are:   $\a1$, $\a3$, $\a4$,
 $\a5$, $\b3$ and $\c1$ with $\c2=\c1$; all  the others vanish.
The equations (\ref{cb})--(\ref{cd}) reduce to a single
relation:  $\a1\a4 +\a5\c1=0$.  The Schouten bracket
$[[r,r]] = -(\a1\b3+\c1^2) K\wedge M\wedge P$ characterizes the
standard and non-standard subfamilies:

\noindent
$\bullet$ Standard subfamily: $\a1$, $\a3$, $\a4$, $\a5$,
$\b3$, $\c1$ with $\c2=\c1$, $\a1\a4 +\a5\c1=0$ and
$\a1\b3+\c1^2\ne 0$.

\noindent
$\bullet$ Non-standard subfamily:   $\a1$, $\a3$, $\a4$,
$\a5$, $\b3$, $\c1$ with $\c2=\c1$, $\a1\a4 +\a5\c1=0$ and
$\a1\b3+\c1^2= 0$.

 The final Lie bialgebra structure and Poisson--Lie brackets for
both subfamilies are given by:
\bea
&&r=\a1 D\wedge P  +\a3 P\wedge M + \a4 H\wedge M +
\a5 P\wedge H  +\b3 K\wedge M  \cr
&&\qquad\qquad  +\c1 ( D\wedge M+P\wedge K)
\nonumber\\[2pt] 
&&\delta(P)=0\qquad \delta(M)=0\label{ne}\\
&& \delta(D)= -\a1 D\wedge P-\a3 P\wedge M - 2 \a4 H\wedge M
- 3 \a5 P\wedge H+\b3 K\wedge M \nonumber \\[2pt]
&&\delta(K)=\a1 (P\wedge K +D\wedge M)  
+\a4 P\wedge M - \a5 H\wedge M  - 2\c1 K\wedge M \nonumber\\[2pt]
&&\delta(H)=
  - 2\a1 P\wedge H  - \b3 P\wedge M + 2 \c1 H\wedge M \nonumber
\\[2pt] &&\delta(C)=   \a1( D\wedge K + 2   P\wedge
C)  + \a3 K\wedge M - \a4 D\wedge M + \a5 (D\wedge P + K\wedge
H)\cr
&&\qquad\qquad - 2 \c1 C\wedge M .
\nonumber
\eea
\bea
&&\{\zd,\zt\}=  0 \qquad
\{\zd,\zx\}=\a1 \left( e^{-\zd} (1-\zc\zt) -1\right)+\a5 \zc
e^{-3\zd} (1-\zc\zt) \nonumber\\[2pt]
&&\{\zd,\zv\}=\a1 \zc e^{-\zd} + \a5 \zc^2 e^{-3\zd} \qquad
\{\zd,\zc\}=0 \nonumber\\[2pt]
&&\{\zt,\zx\}= 2 \a1 \zt + \a5 \left( 1 - 
e^{-3\zd}(1-\zc\zt)\right) \qquad
\{\zt,\zv\}=- \a5 \zc e^{-3\zd} \qquad
\{\zt,\zc\}=0\nonumber\\[2pt]
&&\{\zx,\zv\}=\a1 \zv \qquad
\{\zx,\zc\}=2\a1\zc -  \a5 \zc^2 e^{-3\zd} (1-\zc\zt) \qquad
\{\zv,\zc\}=-\a5  \zc^3 e^{-3\zd}\nonumber\\[2pt]
&&\{\zm,\zd\}=-\a1  \zv  e^{-\zd} (1-\zc\zt)+\a4 
\zc e^{-2\zd} -\a5 \zv  \zc e^{-3\zd} (1-\zc\zt)\nonumber\\[2pt]
 &&\{\zm,\zt\}= \a4 (1-e^{-2\zd})+\a5
\zv e^{-3\zd} (1-\zc\zt)-2 \c1\zt\label{nf}\\
&&\{\zm,\zx\}= \a3 \left(1-e^{-\zd}(1-\zc\zt)
 \right) -\a4 \zv-\frac 12 \a5 \zv^2+\b3 \zt e^{\zd}
\nonumber\\[2pt] &&\{\zm,\zv\}=-\a3 \zc e^{-\zd} -\b3
(e^{\zd}-1)+2\c1 \zv \nonumber\\[2pt] &&\{\zm,\zc\}=\a4 \zc^2
e^{-2\zd} -\a5 \zv \zc^2 e^{-3\zd} (1-\zc\zt) + 2\c1 \zc  .
\nonumber
\eea 

 We recall that the   quantum Schr\"odinger  algebra whose
underlying  non-standard  Lie bialgebra corresponds to setting
$\a1=-z$, $\a3=z/2$ and the remaining parameters equal to zero
was obtained  in \cite{sch1}.
 Furthermore it was proven in \cite{sch3,BHNN} that   a space
discretization of the $(1+1)$ heat-Schr\"odinger equation on a
uniform lattice is endowed with this  Hopf algebra symmetry.


\subsect{H primitive}

  Finally, $\delta(H) = 0$ renders
$\a1$, $\a2$, $\a3$, $\a4$, $\a5$ and $\c2$ with $\b6 =- 2\a1$
 as  free parameters. Moreover the set of equations
(\ref{cb})--(\ref{cd}) implies that $\a1 = 0$ and $\a2\a3 +
\a5\c2 = 0$. The Schouten bracket in this case is given by
$[[r,r]] = - \c2^2 K\wedge M\wedge P$ so that the standard
solution reads

\noindent
$\bullet$ Standard subfamily:  $\a2$, $\a3$, $\a4$, 
$\c2\ne 0$ with $\a5= - {\a2\a3}/{\c2}$.
\bea
&&r=\a2 D\wedge H  +\a3 P\wedge M + \a4 H\wedge M 
 -\frac{\a2\a3}{\c2} P\wedge H   +\c2  P\wedge K
\nonumber\\[2pt] 
&&\delta(H)=0\qquad \delta(M)=0\qquad
\delta(P)=-  \a2 H\wedge P  -  \c2 P\wedge M\label{ng}\\
&& \delta(D)= -2\a2 D\wedge H -\a3 P\wedge M - 2 \a4 H\wedge M
+ 3   \frac{\a2\a3}{\c2} P\wedge H \nonumber \\[2pt]
&&\delta(K)=\a2 (H\wedge K +D\wedge P)  
 +\a4 P\wedge M +\frac{\a2\a3}{\c2} H\wedge M  -  \c2 K\wedge M
\nonumber\\[2pt] 
&&\delta(C)=  2 \a2 H\wedge
C + \a3 K\wedge M - \a4 D\wedge M  -\frac{\a2\a3}{\c2} (D\wedge
P + K\wedge H)  .
\nonumber
\eea
\bea
&&\{\zd,\zt\}=  \a2 ( e^{-2\zd}-1)   \qquad
 \{\zd,\zx\}= \a2\zv  -\frac{\a2\a3}{\c2} \zc
e^{-3\zd} (1-\zc\zt)  \nonumber\\[2pt]
&&\{\zd,\zv\}=   -\frac{\a2\a3}{\c2} \zc^2 e^{-3\zd}  \qquad
 \{\zd,\zc\}=-\a2 \zc^2 e^{-2\zd}\nonumber\\[2pt]
 &&\{\zt,\zx\}=  - \a2 (\zx + 2 \zt\zv)    -\frac{\a2\a3}{\c2}
\left( 1 -  e^{-3\zd}(1-\zc\zt)\right) \nonumber\\[2pt]
&&\{\zt,\zv\}=\a2 \zv+\frac{\a2\a3}{\c2} \zc e^{-3\zd}  \qquad
 \{\zt,\zc\}=2\a2 \zc  \qquad 
 \{\zx,\zv\}= -\a2\zv^2 \nonumber\\[2pt]
&&\{\zx,\zc\}= - 2 \a2 \zv\zc+\frac{\a2\a3}{\c2} \zc^2 e^{-3\zd}
(1-\zc\zt)
\qquad
\{\zv,\zc\}= \frac{\a2\a3}{\c2}  \zc^3 e^{-3\zd} \label{nh}\\
&&\{\zm,\zd\}=-\frac 12 \a2 \zv^2+\a4 
\zc e^{-2\zd} +\frac{\a2\a3}{\c2} \zv  \zc e^{-3\zd}
(1-\zc\zt)\nonumber\\[2pt]
&&\{\zm,\zt\}= \a2 \zt\zv^2+\a4 (1-e^{-2\zd}) 
-\frac{\a2\a3}{\c2}
\zv e^{-3\zd} (1-\zc\zt) \nonumber\\[2pt]
 &&\{\zm,\zx\}= \frac 12 \a2 \zx\zv^2 +\a3
\left(1-e^{-\zd}(1-\zc\zt)
 \right) -\a4 \zv  + \frac{\a2\a3}{2\c2} \zv^2 +\c2
\zx\nonumber\\[2pt] &&\{\zm,\zv\}=-\frac 12 \a2 \zv^3-\a3 \zc
e^{-\zd} +\c2 \zv
\nonumber\\[2pt] 
 &&\{\zm,\zc\}=-\a2 \zc\zv^2+\a4 \zc^2 e^{-2\zd}
+\frac{\a2\a3}{\c2} \zv
\zc^2 e^{-3\zd} (1-\zc\zt).
\nonumber
\eea 

 One arrives at the non-standard subfamily by requiring $\c2=0$,
that is, with parameters $\{\a2,\a3,\a4,\a5\}$, together with 
the constraint
$\a2\a3 = 0$. If we assume $\a2 = 0$, then $\delta(P) = 0$ 
 and we are within the non-standard subfamily of Section 3.2.
Therefore, we have not to take into account this possibility,
and consider  the case
$\a3 = 0$ as follows:

\noindent
$\bullet$ Non-standard subfamily:  $\a2$,   $\a4$, 
$\a5$.  
\bea
&&r=\a2 D\wedge H  + \a4 H\wedge M 
+\a5 P\wedge H   
\nonumber\\[2pt] 
&&\delta(H)=0\qquad \delta(M)=0\qquad 
\delta(P)= -  \a2 H\wedge P \label{ni}\\
&& \delta(D)= -2\a2 D\wedge H  - 2 \a4 H\wedge M
- 3  \a5 P\wedge H \nonumber \\[2pt]
&&\delta(K)=\a2 (H\wedge K +D\wedge P)  
 +\a4 P\wedge M -\a5  H\wedge M 
\nonumber\\[2pt] 
&&\delta(C)=  2 \a2 H\wedge
C  - \a4 D\wedge M +\a5 (D\wedge P +
K\wedge H)  .
\nonumber
\eea
\bea
&&\{\zd,\zt\}=  \a2 ( e^{-2\zd}-1)   \qquad
 \{\zd,\zx\}= \a2\zv +\a5 \zc
e^{-3\zd} (1-\zc\zt)  \nonumber\\[2pt]
&&\{\zd,\zv\}= \a5 \zc^2 e^{-3\zd}  \qquad
 \{\zd,\zc\}=-\a2 \zc^2 e^{-2\zd}\nonumber\\[2pt]
&&\{\zt,\zx\}=  - \a2 (\zx + 2 \zt\zv)  +\a5 \left( 1 - 
e^{-3\zd}(1-\zc\zt)\right) \nonumber\\[2pt]
&&\{\zt,\zv\}=\a2 \zv -\a5 \zc e^{-3\zd}  \qquad
 \{\zt,\zc\}=2\a2 \zc  \qquad 
 \{\zx,\zv\}= -\a2\zv^2 \nonumber\\[2pt]
&&\{\zx,\zc\}= - 2 \a2 \zv\zc-\a5 \zc^2 e^{-3\zd}
(1-\zc\zt)
\qquad
\{\zv,\zc\}= -\a5  \zc^3 e^{-3\zd} \label{nj}\\
&&\{\zm,\zd\}=-\frac 12 \a2 \zv^2+\a4 
\zc e^{-2\zd} -\a5 \zv  \zc e^{-3\zd}
(1-\zc\zt)\nonumber\\[2pt]
 &&\{\zm,\zt\}= \a2 \zt\zv^2+\a4 (1-e^{-2\zd}) +\a5
\zv e^{-3\zd} (1-\zc\zt) \nonumber\\[2pt]
 &&\{\zm,\zx\}= \frac 12 \a2 \zx\zv^2  -\a4 \zv  -\frac 12 \a5
\zv^2
\qquad
\{\zm,\zv\}=-\frac 12 \a2 \zv^3
\nonumber\\[2pt] 
&&\{\zm,\zc\}=-\a2 \zc\zv^2+\a4 \zc^2 e^{-2\zd} -\a5 \zv
\zc^2 e^{-3\zd} (1-\zc\zt).
\nonumber
\eea

 The quantum deformation of the non-standard  Lie bialgebra
corresponding to $\a2=-2z$, $\a4=z$ and $\a5=0$ was constructed
in \cite{sch2}; this
 leads   to a time discretization of the  heat-Schr\"odinger
equation on a uniform lattice with  Hopf algebra symmetry
\cite{sch2,BHNN}.


\sect{Lie sub-bialgebra embeddings}

In this section we investigate  the embeddings (\ref{ab}) at a
Lie bialgebra level, that is,  we analyse which of the
 harmonic oscillator $h_4$ \cite{BH}, extended Galilei $\gal$ 
\cite{gal} and $gl(2)$ \cite{BHP} Lie bialgebras can be embedded
within the   Schr\"odinger  bialgebras as {\em Lie
sub-bialgebras}. In order to state precisely this idea let us
consider  two Lie bialgebras   $(h,\delta_h)$ and 
$(g,\delta_g)$  in such a manner that $h$ is a Lie subalgebra of
$g$: $h\subset g$. We say that $h$ is a Lie sub-bialgebra of $g$
if the cocommutator
$\delta_g$ in $g$ of any generator $X_i\in h$ is of the form
\be
\delta_g(X_i)=f_{i}^{jk} X_j\wedge X_k\qquad X_j,X_k\in h .
\label{jja}
\ee

As a consequence,  the following results establish whether  it
is possible or not to construct a quantum  Schr\"odinger
algebra with one of the above subalgebras as a  Hopf
subalgebra. In other words, whenever the Lie bialgebra
embedding exists for a given subalgebra,  one can make use of 
its known  quantum deformations in order to obtain a  quantum
Schr\"odinger algebra  (this  was the procedure  used for two
particular deformations in
\cite{sch1,sch2}).


\subsect{$h_4$ Lie sub-bialgebras}

 The commutation rules of the harmonic oscillator algebra $h_4$
in the usual basis $\{N,A_+,A_-,M\}$ are given by  
\be
 [N, A_+] = A_+\qquad [N, A_-] = - A_-\qquad [A_-, A_+] = M 
\qquad [M,\,\cdot\,]=0 .
\label{ja}
\ee
 It was proven in \cite{goslar} that all the $h_4$ Lie
bialgebras are coboundary ones; all of them, together
 with their quantum deformations, were explicitly  obtained in
\cite{BH}. The general classical  $r$-matrix for $h_4$, which
depends on six parameters
 $\{\alpha_+,\alpha_-,\beta_+,\beta_-,\vartheta,\xi\}$,  is
given by
\be
 r =\alpha_+ N \wedge A_+ + \alpha_- N \wedge A_- + \vartheta N
\wedge M +
\xi A_+
\wedge
A_- + \beta_+ A_+\wedge M + \beta_- A_- \wedge M
\label{jb}
\ee
where the  parameters must fulfil
\be
 \alpha_+\alpha_- = 0  \qquad \alpha_+(\xi + \vartheta) = 0 
\qquad
\alpha_-(\xi - \vartheta) = 0  
\label{jc}
\ee
and the Schouten bracket reads
\be
 [[r, r]] = ({\alpha}_+{\beta}_- + {\alpha}_-{\beta}_+ -
{\xi}^2) M \wedge A_+ \wedge A_- .
\label{jd}
\ee
 The cocommutators, which can be obtained from (\ref{bb}), turn
out to be
\bea
&&\delta(N)= {\alpha}_+ N \wedge A_+ - {\alpha}_- N \wedge A_-
+ {\beta}_+ A_+ \wedge M - {\beta}_-
A_- \wedge M  \cr
 &&\delta(A_+)= - {\alpha}_-(N \wedge M + A_+ \wedge A_- ) -
(\vartheta +
\xi) A_+\wedge M \cr
 &&\delta(A_-) = {\alpha}_+ ( N \wedge M - A_+ \wedge A_-) +
(\vartheta -
\xi) A_- \wedge M\cr
&&\delta(M) = 0.
\label{je}
\eea

The  $h_4$ algebra can be considered as a subalgebra of
$\Sch$ (\ref{aa})  by denoting the  oscillator generators as
\be
N \to - D   \qquad A_+ \to P   \qquad A_- \to K  \qquad M \to  M .
\label{jf}
\ee
 Now we write (\ref{je}) in terms of $\{D,P,K,M\}$ and according
to (\ref{jja}) impose that their general cocommutators 
(\ref{ci})  give  exactly (\ref{je}). This requirement implies
that
 the Schr\"odinger bialgebra parameters must be
\be
\begin{array}{lll}
  \a1=-{\alpha}_+  & \qquad \a3= \beta_+  &\qquad \c1
=-\vartheta\cr
 \b1=-{\alpha}_-  & \qquad \b3= \beta_-  &\qquad \c2 =\xi
\end{array}
\label{jg}
\ee
with  all the remaining ones  equal to zero.
 Under these assumptions  it can be checked that the full  set
 of 19 equations (\ref{cb})--(\ref{cd}) reduces to (\ref{jc}),
the Schr\"odinger $r$-matrix (\ref{ca}) is  equal to (\ref{jb})
and the Schouten bracket (\ref{ce}) reproduces (\ref{jd}). Hence
we conclude that

\noindent
 {\bf Proposition 4.1.} {\sl   All harmonic oscillator Lie
bialgebras with generators $\{D,P,K,\break  M\}$ are
Schr\"odinger Lie sub-bialgebras. The resulting Schr\"odinger
classical $r$-matrix and cocommutators  containing the $h_4$
sub-bialgebras  depend on six parameters
 $\{\alpha_+,\alpha_-,\beta_+,\beta_-,\vartheta,\xi\}$ 
subjected to  (\ref{jc}); they are
\be 
\begin{array}{l}
r=-\alpha_+ D\wedge P  
 - \alpha_- D\wedge K   +
\beta_+ P\wedge M   +
\beta_- K\wedge M    -\vartheta D\wedge M + \xi P\wedge K \\[6pt]
\delta(D)= {\alpha}_+ D \wedge P - {\alpha}_- D \wedge K
- {\beta}_+ P \wedge M + {\beta}_-
K \wedge M  \cr
 \delta(P)=   {\alpha}_-(D \wedge M -P \wedge K ) - (\vartheta +
\xi) P\wedge M \cr
 \delta(K) = -{\alpha}_+ ( D \wedge M + P \wedge K) + (\vartheta
- \xi) K \wedge M\cr
\delta(M) = 0\cr
\delta(H)=
   2\alpha_+ P\wedge H 
+ \alpha_- ( D\wedge P
+ 2  K\wedge H ) - \beta_- P\wedge M 
- 2 \vartheta H\wedge M \cr
 \delta(C)= - 2\alpha_- K\wedge C-  \alpha_+ ( D\wedge K+ 2  
P\wedge C)
 + \beta_+ K\wedge M  + 2 \vartheta C\wedge M .
\end{array}
\label{jjg}
\ee
 The Schouten bracket is 
 $[[r, r]] = ({\alpha}_+{\beta}_- + {\alpha}_-{\beta}_+  -
{\xi}^2) K
\wedge M \wedge P$.}


\subsect{$gl(2)$ Lie sub-bialgebras}

Now, we consider the $gl(2)$ Lie algebra whose  generators
$\{J_3,J_+,J_-,I\}$ satisfy the following  commutation rules:
\be
 [J_3, J_+]= 2 J_+  \qquad [J_3, J_-] = - 2 J_-  \qquad [J_+, 
J_-] = J_3 
\qquad [I,\,\cdot\,]=0 
\label{jh}
\ee
 where $I$ is the central generator. All the $gl(2)$ Lie
bialgebras are coboundary ones  coming from a classical
$r$-matrix which depends on six parameters
$\{a_+,a_-,b_+,b_-,a,b\}$ \cite{BHP}:
\be
 r=\frac 12(a_+ J_3 \wedge J_+ - a_- J_3 \wedge  J_- - b J_3
\wedge I + b_+ J_+ \wedge I - b_-  J_- \wedge I - 2a J_+ \wedge 
J_-)
\label{ji}
\ee
which is subjected to the relations
\be
a_+b - b_+a = 0\qquad a_+b_- + a_-b_+ = 0  \qquad a_-b +
b_-a = 0. 
\label{jj}
\ee
The Schouten bracket is given by
\be
[[r, r]] = (a^2 + a_+a_-)  J_3 \wedge J_+ \wedge J_- .
\label{jk}
\ee
 The corresponding cocommutators are  
\bea
 &&\delta(J_3)= a_+ J_3 \wedge J_+ + a_- J_3 \wedge  J_- + b_+
J_+
\wedge I + b_-  J_- \wedge I \cr
 &&\delta(J_+) = a J_3 \wedge J_+ - \frac{1}2 b_- J_3 \wedge I +
a_- J_+ \wedge  J_- + b J_+ \wedge I \cr
 &&\delta( J_-) = a J_3 \wedge  J_- - \frac{1}2 b_+ J_3 \wedge I
- a_+ J_+ \wedge  J_- - b  J_- \wedge I \cr
&&\delta(I) = 0.
\label{jl}
\eea

The $gl(2)$ algebra can be regarded as a subalgebra of $\Sch$
 once we   rename the $gl(2)$ generators as:
\be
 J_3\to -D \qquad J_+ \to H \qquad   J_-  \to -C \qquad    I\to
M  .
\label{jm}
\ee
We demand now that the Schr\"odinger cocommutators
(\ref{ci})  for $\{D,H,C,M\}$ lead to (\ref{jl}); thus  the
Schr\"odinger bialgebra parameters turn out to be 
\be
\begin{array}{lll}
 \a2=-\frac 12 a_+  & \qquad \a4=\frac 12 b_+  &\qquad \c1
=\frac 12 b\\[4pt]
 \b2=-\frac 12 a_-  & \qquad \b4= \frac 12 b_-  &\qquad \c3 =a
\end{array}
\label{jn}
\ee
with $\c2$ arbitrary and all the others equal to zero. 
In this case,   the  equations
(\ref{cb})--(\ref{cd}) give rise to (\ref{jj}), together with
\be
a^2 + a_+a_-=0.
\label{jo}
\ee
  This last equation implies that the  Schouten  bracket of
$gl(2)$  vanishes, which in turn precludes the embedding of the
standard $gl(2)$ bialgebras. These results can be summarized as
follows:

\noindent
 {\bf Proposition 4.2.} {\sl Amongst all the $gl(2)$ Lie
bialgebras with generators
 $\{D,H,C, M\}$, only the  non-standard  ones are Schr\"odinger
Lie sub-bialgebras. The Schr\"odinger $r$-matrix and
cocommutators which comprise the  non-standard $gl(2)$
sub-bialgebras depend on seven parameters
 $\{a_+,a_-,b_+,b_-,a,b,\c2\}$ which fulfil (\ref{jj}) and
(\ref{jo}):
\be 
\begin{array}{l}
r=\frac 12\bigl( -a_+ D\wedge H 
 - a_- D\wedge C   +
b_+ H\wedge M   +
b_- C\wedge M   \cr
\qquad\qquad\qquad  +b D\wedge M +2 a H\wedge C +2 \c2 P\wedge K
\bigr)\\[6pt]
\delta(D)= a_+ D \wedge H - a_- D \wedge C
- b_+ H \wedge M + b_-
C \wedge M  \cr
 \delta(H)= - a   D \wedge H + \frac 12 b_-   D\wedge M - a_-
H\wedge C + b H \wedge M\cr
 \delta(C) = - a   D \wedge C - \frac 12 b_+   D\wedge M - a_+
H\wedge C - b C \wedge M\cr
\delta(M) = 0\cr
 \delta(P)=\frac 12\bigr(a_+ H\wedge P + a_- (C\wedge P +
D\wedge K)  - b_- K\wedge M \cr
\qquad\qquad\qquad  +(b-2\c2) P\wedge M + 2a K\wedge H\bigl)\cr
 \delta(K)=-\frac 12\bigr(a_- C\wedge K + a_+ (H\wedge K +
D\wedge P) - b_+ P\wedge M
\cr
\qquad\qquad\qquad  +(b+2\c2) K\wedge M -2 a P\wedge C \bigl).
\end{array}
\label{jp}
\ee
The Schouten bracket is 
$[[r, r]] =  - \c2^2 K \wedge
M \wedge P$.}

This statement shows that there exist both 
 standard Schr\"odinger Lie bialgebras ($\c2\ne 0$) and
non-standard ones ($\c2=0$) that contain non-standard 
 $gl(2)$ Lie sub-bialgebras. On the other hand, 
let us recall that
 in \cite{Dobrev}, a $q$-deformed Schr\"odinger algebra was
obtained by demanding that the subalgebra structure of
(standard)  $sl_q(2)$ should be preserved by the deformation.
This led to a deformation of
$\Sch$ for which   the previous analysis  proves that no  Hopf
structure can be constructed.  In other words, there cannot
exist any Hopf algebra deformation of
$\Sch$ retaining the standard $sl_q(2)$ as a Hopf subalgebra.
 This example enlightens the usefulness of a Lie bialgebra
approach when looking for a precise quantum deformation of a
given Lie algebra.


\subsect{$\gal$ Lie sub-bialgebras}

In the usual kinematical basis, the $(1+1)$-dimensional extended
Galilei algebra $\gal$  is spanned by $\{K,H,P,M\}$,  directly
arising  as a subalgebra of $\Sch$ (\ref{aa}).  The
cocommutators of the $\gal$ Lie bialgebras can be written
collectively in terms of  nine parameters
$\{\alpha,\xi,\nu,\beta_1,
\beta_2,\beta_3,\beta_4,\beta_5,\beta_6\}$ as
\cite{gal} (see also \cite{anna}):
\bea
&&\delta(K)=\beta_6 K\wedge P+ \xi K\wedge M + \nu P\wedge H +
\beta_1 P\wedge M + \beta_2 H\wedge M\cr
&&\delta(H)=\beta_5 K\wedge M -(\beta_6+\alpha)  P\wedge H +
\beta_3 P\wedge M + (\beta_4 - \xi) H\wedge M\cr
&&\delta(P)=\beta_4 P\wedge M + (\beta_6 + \alpha) H\wedge M\cr
&&\delta(M)=\alpha P\wedge M  
\label{jq}
\eea
 satisfying the  following equations:
\bea
&&\alpha\beta_5=0\qquad
\beta_6(\beta_6 +\alpha)=0\qquad
\beta_4(\beta_6 +\alpha)=0\cr
&&\nu(\xi -\beta_4)=0\qquad
\alpha(\xi -\beta_4) - \nu\beta_5=0 .
\label{jr}
\eea
 Unlike the $h_4$ and $gl(2)$ cases, the $\gal$ Lie bialgebras
include both non-coboundary and  coboundary Lie bialgebras 
\cite{gal}. The standard $\gal$ bialgebras depend on two
 parameters $\{\xi\ne 0,\beta_1\}$ with $\beta_4=\xi$ while the
others are equal to zero; the corresponding classical $r$-matrix
reads
\be
r=\xi K\wedge P + \beta_1 H\wedge M .
\label{jjr}
\ee
The non-standard ones correspond to take the three parameters
$\{\beta_1,\beta_2,\beta_3\}$ and the $r$-matrix is given by
\be
r=  \beta_1 H\wedge M+ \beta_2 H\wedge P+ \beta_3 M\wedge K .
\label{jjs}
\ee

We impose the 
Schr\"odinger cocommutators
(\ref{ci}) to be equal to (\ref{jq}); this implies that
 $\alpha$, $\nu$, $\beta_5$ must vanish, while the Schr\"odinger
parameters are
\be
\begin{array}{lll}
  \a4=\beta_1 & \qquad \a5=-\beta_2  &\qquad \c1=\frac 12
(\beta_4 -\xi)\\[4pt]
  \b3=-\beta_3  & \qquad \b6=\beta_6  &\qquad \c2=-\frac 12
(\beta_4 +\xi)
\end{array}
\label{js}
\ee
 with $\a3$ arbitrary and the remaining ones equal to zero. The
equations  (\ref{cb})--(\ref{cd}) lead to $\beta_6=0$ and
\be
\beta_2 (2\beta_4 -\xi)=0 .
\label{jt}
\ee
 Therefore the equations (\ref{jr}) are trivially satisfied. 
The final result is summed up by

\noindent
 {\bf Proposition 4.3.} {\sl Amongst all the extended Galilei 
Lie bialgebras with generators $\{K,H,P, M\}$ and cocommutators
given by (\ref{jq}), only those with
$\alpha=\nu=\beta_5=\beta_6=0$ are Schr\"odinger Lie
sub-bialgebras. The
 Schr\"odinger Lie  bialgebras which   contain the $\gal$
sub-bialgebras depend on six parameters
 $\{\xi,\beta_1,\beta_2,\beta_3,\beta_4,\a3\}$  fulfilling
(\ref{jt}):
\be 
\begin{array}{l}
 r=  \a3 P \wedge M +\beta_1 H \wedge M - \beta_2 P \wedge H -
\beta_3 K
\wedge M\cr
\qquad\qquad\qquad  + \frac 12 (\beta_4 - \xi) D \wedge M -
\frac 12 (\beta_4 + \xi) P \wedge K \\[6pt]
\delta(K)= \xi K\wedge M   +
\beta_1 P\wedge M + \beta_2 H\wedge M\cr
\delta(H)=  
\beta_3 P\wedge M + (\beta_4 - \xi) H\wedge M\cr
\delta(P)=\beta_4 P\wedge M  \cr
\delta(M)=0\cr
 \delta(D)=- \a3 P \wedge M - 2 \beta_1 H \wedge M + 3 \beta_2 P
\wedge H -\beta_3 K \wedge M\cr
 \delta(C) = \a3 K \wedge M - \beta_1 D \wedge M -\beta_2 (D
\wedge P + K \wedge H)-(\beta_4 -\xi) C \wedge M .
\end{array}
\label{ju}
\ee
The Schouten bracket is 
$[[r, r]] =  -\frac 14 (\beta_4 +\xi )^2 K \wedge
M \wedge P$.}

 Notice that {\em all} the coboundary $\gal$ bialgebras
determined by either (\ref{jjr}) or (\ref{jjs}) are 
Schr\"odinger   sub-bialgebras. There are also non-coboundary
$\gal$ bialgebras arising as Schr\"odinger  sub-bialgebras.


 \sect{Quantum Schr\"odinger algebras and quantum universal
$R$-matrices}

 The main advantage of a systematic analysis of  Lie bialgebra
structures on a given Lie algebra is that they  characterize 
possible `directions' for the obtention of  quantum deformations
which  are guided by the cocommutator $\delta$, that gives the
first order in the deformation of the comultiplication map. As
we have already mentioned, there are only two  quantum
$\Sch$ algebras endowed with a known Hopf structure
 \cite{sch1,sch2}, and both of them are of non-standard type. 
In this section we focus mainly on the standard $\Sch$
bialgebras with a two-fold purpose. First, we construct  new
quantum Schr\"odinger algebras starting from Lie bialgebras with
either $D$ or $H$ as primitive generators. Second, we illustrate
with these examples the Lie sub-bialgebra embeddings studied in
the previous section and    exploit them in order to deduce 
universal
$R$-matrices for these quantum $\Sch$ algebras.


\subsect{D primitive: $U_{\c1,\c2}(\Sch)$}

Let us consider the  Schr\"odinger bialgebras with $D$ primitive
that are defined through (\ref{na}). The  corresponding 
two-parameter  Hopf algebra $U_{\c1,\c2}(\Sch)$ has the
following coproduct and  commutation rules:
\bea
 &&\Delta(D) = 1\otimes D + D \otimes 1  \qquad \Delta(M) = 1
\otimes M + M
\otimes 1  \cr
&&\Delta(P) = 1 \otimes P + P \otimes e^{(\c1 - \c2) M}\qquad
\Delta(K) = 1 \otimes K + K \otimes e^{-(\c1 + \c2) M}\cr
&&\Delta(H) = 1 \otimes H + H \otimes e^{2\c1 M}\qquad
\Delta(C) = 1 \otimes C + C \otimes e^{-2\c1 M}  
\label{ka}
\eea
\be
\begin{array}{lll}
 [D, P]=-P \quad  &[D, K]=K\quad &\displaystyle{[K,
P]=\frac{1-e^{-2\c2 M}}{2\c2}}  \cr 
[D, H]=-2H \quad  &[D, C]=2C \quad &[H, C]= D\cr
[K, H]=P \quad  &[K, C] =0 \quad
&[M,\,\cdot\,]=0\cr
[P, C]=-K\quad  &[P, H]=0 .\quad
& \cr
\end{array}
\label{kb}
\ee
The counit is  trivial and the antipode can be easily deduced
from the Hopf algebra axioms. Notice that $U_{\c1,\c2}(\Sch)$
is a standard quantum algebra  whenever $\c2\ne 0$, while it
leads to a non-standard one under the limit $\c2\to 0$, 
$U_{\c1}(\Sch)$, which has non-deformed commutation rules.

By taking into account the results presented in  section 4, it
can be checked that the Lie bialgebra (\ref{na}) has $h_4$,
$gl(2)$ and $\gal$ Lie sub-bialgebras. After deformation, these
structures  give rise to  Hopf subalgebras of
$U_{\c1,\c2}(\Sch)$, namely

\noindent
 (i) A standard quantum algebra $U_{\vartheta,\xi}(h_4)$ of type
II of \cite{BH}, with generators $\{D,P,\\ K,M\}$ and  
$\vartheta = -\c1$, $\xi =  \c2$.

\noindent
 (ii) A non-standard quantum algebra  $U_b(gl(2))$ belonging to
the family II of \cite{BHP}, with generators $\{D,H,C, M\}$ and 
$b = 2\c1$.

\noindent
 (iii) A non-coboundary  quantum algebra $U_{\xi,\beta_4}(\gal)$
of the family Ia of \cite{gal}, with generators $\{K,H,P, M\}$ 
and
$\xi=-\c1-\c2$,  $\beta_4=\c1-\c2$.

 In what follows we show how the  above quantum subalgebras
allow us to obtain directly universal $R$-matrices for either
$U_{\c2}(\Sch)$ or
$U_{\c1}(\Sch)$.

  We consider the standard  algebra  $U_{\c2}(\Sch)$ by fixing  
  $\c1=0$  in (\ref{ka}) and (\ref{kb}). It contains a standard
quantum extended Galilei subalgebra $U_{\xi}(\gal)$ (note that
in this case
$\beta_4=\xi=-\c2$) whose quantum universal $R$-matrix was
deduced in
\cite{gal}. In the Schr\"odinger basis this element reads
\bea
&&{\cal R}=\exp\{\c2 P\wedge K f(M,\c2)\}\nonumber\\[2pt]
&&f(M,\c2)=\frac{e^{\c2 M/2}\otimes e^{\c2 M/2}}{\sqrt{\sinh
\c2 M\otimes \sinh \c2 M}}\arcsin\left({\frac{\sqrt{\sinh \c2
M\otimes
\sinh \c2 M}}{\cosh((\c2/2)\Delta(M))}}\right) .
\label{kc}
\eea
As it was proven in \cite{gal},  ${\cal R}$  is
not a solution of the quantum YBE but it fulfils 
\be
{\cal R}\Delta(X){\cal R}^{-1}=\sigma\circ\Delta(X) ,
\label{kd}
\ee
 where  $\sigma(X\otimes Y)=Y\otimes X$, for the Galilei 
generators
$\{K,H,P, M\}$. Furthermore, since
\be
[\c2 P\wedge K f(M,\c2),1\otimes X + X\otimes 1]=0\quad
{\mbox {for}}\quad X=\{D,C\}
\ee
the relation (\ref{kd})  is also satisfied by  $D$
 and $C$   so that (\ref{kc})  is a non-quasitriangular 
universal $R$-matrix for $U_{\c2}(\Sch)$. Notice that its
underlying classical
$r$-matrix is $r=\c2 P\wedge K$.

 On the other hand, the non-standard quantum algebra 
$U_{\c1}(\Sch)$ provided by   $\c2\to 0$  comprises a
non-standard quantum subalgebra  $U_b(gl(2))$ whose universal
$R$-matrix was obtained in \cite{BHP}:
\be
{\cal R}=\exp\{-\c1 M\otimes D \}\exp\{\c1 D\otimes M \} 
\label{ke}
\ee
which not only verifies the relation (\ref{kd}) but also
the quantum YBE   for
$\{D,H,C, M\}$. The two remaining generators
$P$ and $K$ also fulfil (\ref{kd}):
\bea
 &&\exp\{\c1 D\otimes M \}\Delta(X)\exp\{-\c1 D\otimes M
\}=1\otimes X+ X\otimes 1\equiv \Delta_0(X)\cr
 &&\exp\{-\c1 M\otimes D \}\Delta_0(X)\exp\{\c1 M\otimes D
\}=\sigma\circ 
\Delta(X)\ \ \mbox{for}\ X=\{P,K\} .
\label{kke}
\eea
Consequently, the element (\ref{ke}) is a
triangular universal $R$-matrix for $U_{\c1}(\Sch)$ with classical
$r$-matrix given by
$r=\c1 D\wedge M$.


\subsect{H primitive: $U_{\a2,\c2}(\Sch)$}

 As a second example, we consider the standard $\Sch$ bialgebra
with $H$ primitive (\ref{ng}) with two  parameters $\a2$,
$\c2\ne 0$ while $\a3$,
$\a4$ are taken to be zero:
\bea
&&r =\a2 D \wedge H + \c2 P \wedge K\cr
&&\delta(H) = 0 \qquad \delta(M) = 0\cr
&&\delta(D) = - 2\a2 D \wedge H   \qquad 
\delta(P) =  P \wedge (\a2 H - \c2  M)  \cr
&&
\delta(C) = - 2\a2 C \wedge H  \qquad
\delta(K) =   K \wedge (- \a2 H  - \c2   M)+  \a2 D \wedge P .
\label{kf}
\eea
The quantum deformation of this bialgebra  leads  to a
two-parameter Hopf algebra $U_{\a2,\c2}(\Sch)$ endowed with
the following  coproduct and commutation rules:
\bea
&&\Delta(H) = 1 \otimes H + H \otimes 1\qquad
\Delta(M) = 1 \otimes M + M \otimes 1\cr
&&\Delta(D) = 1\otimes D + D \otimes e^{- 2\a2 H}  \qquad 
\Delta(C) = 1 \otimes C + C \otimes e^{- 2\a2 H}  \cr
&&\Delta(P) = 1 \otimes P + P \otimes e^{\a2 H} e^{ - \c2 M}\cr
&&\Delta(K) = 1 \otimes K + K \otimes e^{- \a2 H} e^{-\c2 M}
+\a2 D \otimes e^{- 2\a2 H} P
\label{kg}
\eea
\bea
 &&[D, P] = - P   \qquad [D, K] = K  \qquad [K, P] =  \frac{1 -
e^{- 2\c2 M}} {2\c2}  \cr
 &&[D, H] = \frac{e^{- 2\a2 H} - 1}{\a2}  \qquad [D, C] = 2C -
\a2 D^2 
\qquad [H, C] = D  \cr
 &&[K, H] = e^{ - 2\a2 H} P  \qquad [K, C] = - \frac{1}2 \a2(K 
D+ DK)  
\qquad [M,\,\cdot\,]=0\cr 
&&[P, C] = - K + \frac{1}2 {\a2}  (D P + P D)     \qquad
[P, H] = 0. 
\label{kh}
\eea
 This quantum algebra is of  standard type whenever $\c2\ne 0$;
however if we perfom the limit $\c2\to 0$ in (\ref{kg}) and
(\ref{kh}),  then we obtain a non-standard quantum algebra
$U_{\a2}(\Sch)$ whose Lie bialgebra belongs to (\ref{ni}).

The generators $\{D,H,C,M\}$  close a non-standard 
Hopf subalgebra
$$
 U_{a_+}(gl(2))=U_{a_+}(sl(2,\R))\oplus u(1)\subset
U_{\a2,\c2}(\Sch)
$$
belonging to  the family $\mbox{I}_+$ of \cite{BHP}, where
 $a_+=-2\a2$, $u(1)$ corresponds to the central generator $M$,
and $U_{a_+}(sl(2,\R))$ is the well-known Jordanian or
non-standard quantum $sl(2,\R)$ algebra
\cite{Demidov,Zak,nonsc,Aizawa,vander} writen in the basis of
\cite{nonsf}.  The quantum universal $R$-matrix of
$U_{a_+}(sl(2,\R))$ \cite{nonsf} (see also \cite{nonse}) can be
expressed in the Schr\"odinger basis as
\be
{\cal R}=\exp\{-\a2 H\otimes D \}\exp\{\a2 D\otimes H\} .
\label{ki}
\ee
In this respect, we stress that the  quantum universal
$R$-matrix for $U_{a_+}(sl(2,\R))$ was firstly obtained by
Ogievetsky through a twist operator  \cite{Ogia}.
The element ${\cal R}$ is a solution of the quantum YBE and
also fulfils (\ref{kd}) for $\{D,H,C,M\}$. If we restrict now
to the non-standard quantum algebra $U_{\a2}(\Sch)$ by taking
$\c2\to 0$, it  can be checked that the generators $P$ and $K$
also satisfy the relation (\ref{kd}) by use of similar equations
to (\ref{kke}). Hence we conclude that (\ref{ki}) is a
triangular universal $R$-matrix for $U_{\a2}(\Sch)$ whose
corresponding   classical  $r$-matrix is $r=\a2 D\wedge H$.


\sect{Concluding remarks}

 This paper presents the explicit computation  and analysis of
all the possible Lie bialgebra structures associated to the
Schr\"odinger algebra, a non-semisimple Lie algebra of dimension
6. In general, this kind of classification is shown to be
feasible and enables a systematic study of all possible quantum 
deformations of a given Lie algebra. In particular, this kind of
construction simplifies quite efficiently the task of finding a
certain deformation that should fit with some fixed properties
that are known {\it a priori}. 

 In particular, the coboundary nature of all  the Schr\"odinger
Lie bialgebras is a remarkable result that certainly enhances
the structural vicinity of this Lie algebra with respect to
semisimple ones. On the other hand, the richness of the
subalgebra structure of $\Sch$ allows to show a great variety of
situations concerning the generalization of Lie sub-bialgebras
to full Lie bialgebras on $\Sch$. In this respect, we find all
kinds of situations: sometimes, standard Lie sub-bialgebras can
be generalized, and sometimes they cannot; moreover, even
non-coboundary Lie sub-bialgebras can be generalized to
coboundary ones in the bigger Lie algebra. These examples
suggest that the embedding of a given non-semisimple algebra
within a higher dimensional one is an interesting procedure that
could provide some new information connecting quantum
deformations of both algebras.

On the other hand, we would like to recall that the
two-photon algebra $h_6$ \cite{Gil}, generated by the
 operators $\{\aaa,\ap,\am,\bp,\bm,\bb\}$ and endowed with the
following commutation rules
\be
\begin{array}{lll}
 [\aaa,\ap]=\ap \quad  &[\aaa,\am]=-\am \quad &[\am,\ap]=\bb \cr
 [\aaa,\bp]=2\bp \quad  &[\aaa,\bm]=-2\bm \quad
&[\bm,\bp]=4\aaa+2\bb
\cr
 [\ap,\bm]=-2\am \quad  &[\ap,\bp]=0 \quad
&[\bb,\,\cdot\,]=0\cr
 [\am,\bp]=2\ap \quad  &[\am,\bm]=0 \quad
& \cr
\end{array}
\label{aaz}
\ee
is isomorphic to $\Sch$ through the map
\be
D=-\aaa-\frac 12\bb \qquad P=\ap\qquad K=\am\qquad
H=\frac 12\bp\qquad C=\frac 12\bm .
\label{fb}
\ee
Therefore, all the results contained  in this paper can be
immediately translated into the two-photon algebra language 
through (\ref{fb}), and the new quantum deformations that have
been constructed could be interpreted in a completely
different physical framework. Recall that the two-photon
algebra is the dynamical symmetry algebra of the single-mode
radiation field Hamiltonian that describes in a unified setting
coherent, squeezed and intelligent states of light \cite{Brif}.
Thus, the results here presented provide a
preliminary basis for the application of quantum
deformations in the construction of non-classical states of
light (see, for instance, 
\cite{sch1,sch2} where deformed Fock--Bargmann realizations
corresponding to two different quantum deformations of the
two-photon/Schr\"odinger algebra were introduced).

Finally,  we would like to emphasize
that quantum algebras and  
completely integrable systems are directly related through the
formalism introduced in \cite{BR}. In particular, the
Hamiltonian 
\be
{\cal H}_{gl(2)}=\sum_{i=1}^N\frac {p_i^2}{2\masa_i}+{\cal
G}\left(\sum_{i=1}^N \masa_i q_i^2\right),
\label{finc}
\ee
where ${\cal G}$ is an arbitrary function, can be constructed
through the Poisson realization of the
$gl(2)=\{N,B_+,B_-,M\}$ coalgebra, and 
completely
integrable (and  superintegrable) deformations of
(\ref{finc}) were given through some quantum
deformations of $gl(2)$ in \cite{chains}. On the other hand,
completely integrable deformations of  Hamiltonians of the type
\be
{\cal H}_{\gal}=\sum_{i=1}^N\frac {p_i^2}{2\masa_i}+{\cal
F}\left(\sum_{i=1}^N\masa_i q_i\right) ,
\label{finb}
\ee
where ${\cal F}$ is again arbitrary, can be obtained
through  quantum Poisson  $\gal$ coalgebras (with
generators 
$\gal=\{B_+,A_+,A_-,M\}$) \cite{gal}. 

Since both $gl(2)$ and $\gal$ are subalgebras of $\Sch$, these
constructions can be generalized by making use of
Schr\"odinger coalgebras. In fact, by following
\cite{BR}, it can be shown that the non-deformed
Schr\"odinger coalgebra gives rise to the $N$-particle
Hamiltonian
\be
{\cal H}_{\Sch}=\sum_{i=1}^N\frac {p_i^2}{2\masa_i}+{\cal
F}\left(\sum_{i=1}^N\masa_i q_i\right)
+{\cal
G}\left(\sum_{i=1}^N \masa_i q_i^2\right) 
\label{find}
\ee
through the following phase space
realization of $\Sch$
\be 
 N=q p-\frac {\masa}{2}\qquad A_+=p\qquad A_-=\masa q\qquad
 M=\masa\qquad B_+=\frac{p^2}{\masa}\qquad B_-=\masa q^2 .
\label{fina}
\ee 
Integrals of the motion for (\ref{find}) can be
obtained by using the coproducts of the (fourth-order)
Casimir of
$\Sch$
\cite{Patera}. From this perspective, each
quantum deformation of the Schr\"odinger algebra provides an
integrable deformation of the Hamiltonian (\ref{find}). In
order to obtain such deformations explicitly, the
corresponding deformed Casimir and phase space realization
is needed. A comprehensive treatment of the
integrability properties linked to (classical and quantum)
Schr\"odinger coalgebras will be given elsewhere.


\section*{Acknowledgments}
\noindent
 A.B. and
F.J.H. have been partially supported by  Junta
de Castilla y Le\'on, Spain  (Project   CO2/399).
 P.P. has been supported by a fellowship from AECI, Spain.

\bigskip
\bigskip


\noindent
{\Large{{\bf Appendix: proof of Theorem 2.1}}}

\appendix

\setcounter{equation}{0}

\renewcommand{\theequation}{A.\arabic{equation}}

\bigskip

\noindent
 As $\Sch$ is a six-dimensional Lie algebra, the most general
cocommutator
$\delta: \Sch\rightarrow \Sch \otimes \Sch$ will
be a linear combination (with $6 \times 15 $ real coefficients)
of skewsymmetric products of  the generators $X_i$ of  $\Sch$:
\be
 \delta(X_i)=f_i^{jk}\,X_j\wedge X_k \qquad j < k\qquad
i,j,k=1,\dots,6 
\label{apenda}
\ee
which in matrix form reads
\be
\delta\left(\begin{array}{c}
X_1 \cr X_2 \cr X_3 \cr X_4 \cr X_5 \cr X_6
\end{array}\right)=
\left(\begin{array}{cccccc}
f_1^{12} &\dots & f_1^{16} & f_1^{23} &\dots& f_1^{56} \cr
f_2^{12} & \dots & f_2^{16} & f_2^{23} & \dots& f_2^{56}  \cr
f_3^{12} & \dots & f_3^{16} & f_3^{23} &\dots& f_3^{56}  \cr
f_4^{12} & \dots & f_4^{16} & f_4^{23} & \dots& f_4^{56}  \cr
f_5^{12} &\dots & f_5^{16} & f_5^{23} & \dots& f_5^{56}  \cr
f_6^{12} & \dots & f_6^{16} & f_6^{23} &\dots& f_6^{56}
\end{array}\right)
{\left(\begin{array}{c}
X_1 \wedge X_2  \cr \vdots\cr X_1 \wedge X_6 \cr X_2 \wedge
X_3 \cr \vdots\cr X_5 \wedge X_6
\end{array}\right)} .
\label{apendb}
\ee
We choose 
\be
 X_1=D\qquad X_2=C\qquad X_3=H\qquad X_4=K\qquad X_5=P\qquad
X_6=M
\ee
 and impose the  cocycle condition (\ref{ba}) onto the arbitrary
expression (\ref{apendb}). This implies that the initial 90
coefficients $f_i^{jk}$  can be expressed in terms of 15
independent parameters that we denote
$\aa{l}$
$(l=1,\dots,15)$, namely
\bea
&&f_1^{12} = \aa1  \qquad f_1^{13} = \aa2  \qquad f_1^{14}
= \aa3  \qquad f_1^{15} = \aa4\qquad f_2^{12} = \aa5 \cr
&&f_2^{26} = \aa6\qquad f_1^{24} = \aa7\qquad 
 f_1^{25} = \aa8  \qquad f_1^{26} = \aa9  \qquad f_1^{34} =
\aa{10}\cr
 &&f_1^{35} = \aa{11} \qquad  f_1^{36} = \aa{12}\qquad 
 f_4^{46} = \aa{13}  \qquad f_1^{46} = \aa{14}  \qquad f_1^{56}
= \aa{15} .
\label{apendc}
\eea
The remaining  non-vanishing coefficients read
\bea
&& f_2^{14} = \aa{10} - \aa4\qquad
 f_2^{15} =\aa{11}/3  \qquad f_2^{16} =  \aa{12}/2\qquad
f_2^{23} = \aa2  \cr
&& f_2^{24} = \aa8 - 2\aa3\qquad
f_2^{25} = 2\aa4  \qquad 
f_2^{34} = -\aa{11}/3  \qquad f_2^{46} = -\aa{15}  \cr
&&f_3^{13} = \aa5\qquad
f_3^{14} = \aa7/3\qquad
f_3^{15} = \aa8 - \aa3  \qquad f_3^{16}= \aa9/2  \cr
&& f_3^{23} =
-\aa1\qquad
 f_3^{25} = -\aa7/3\qquad
f_3^{34} = 2\aa3  \qquad f_3^{35} = \aa{10} - 2\aa4  \cr
&& f_3^{36} =
-\aa6\qquad
   f_3^{56} = -\aa{14}\qquad
 f_4^{15} = -\aa2/2  \qquad f_4^{16} = -\aa4  \qquad f_4^{24} =
\aa1/2 
\cr
&&f_4^{25} = \aa5  \qquad
f_4^{26} = \aa8  \qquad f_4^{34} = -\aa2/2  \qquad f_4^{36} =
-\aa{11}/3  \cr
&& f_4^{45} = \aa{10} + \aa4 \qquad f_4^{56} = -\aa{12}/2 
 \qquad f_5^{14} =
-\aa1/2  \qquad f_5^{16} = -\aa3  \cr
 &&f_5^{25} = -\aa1/2  \qquad f_5^{26} = -\aa7/3  \qquad
f_5^{34} = \aa5 
\qquad f_5^{35} = \aa2/2  \cr
 &&f_5^{36} = \aa{10}  \qquad f_5^{45} = -(\aa3 + \aa8)  \qquad
f_5^{46} = -\aa9/2  \qquad f_5^{56} = \aa{13} - \aa6.
\label{apendd}
\eea
Hence we have obtained the following  cocommutator:
\bea
&&\delta(D)= \aa1 D\wedge C + \aa2 D\wedge H + \aa3 D \wedge K +
\aa4 D \wedge P \cr
&&\qquad\quad + \aa7 C \wedge K
+ \aa8 C\wedge P + \aa9 C\wedge M + \aa{10} H\wedge K \cr
&&\qquad\quad+ \aa{11} H\wedge P + \aa{12} H\wedge M
+\aa{14} K\wedge M + \aa{15} P\wedge M \cr
&& \delta(C)= \aa5 D\wedge C + (\aa{10} - \aa4) D \wedge K +
\frac{\aa{11}}3 D\wedge P \cr
&&\qquad\quad+ \frac{\aa{12}}2 D\wedge M
+ \aa2 C\wedge H + (\aa8 - 2\aa3) C \wedge K
+ 2\aa4 C\wedge P\cr
&&\qquad\quad + \aa6 C\wedge M
-\frac{\aa{11}}3 H\wedge K
- \aa{15} K\wedge M  \cr
&& \delta(H)=\aa5 D\wedge H + \frac{\aa7}3 D \wedge K +
(\aa8 - \aa3) D \wedge P \cr
&&\qquad\quad+ \frac{\aa9}2 D\wedge M
-\aa1 C\wedge H -
\frac{\aa7}3 C\wedge P + 2\aa3 H\wedge K \cr
&&\qquad\quad+ (\aa{10} - 2 \aa4) H\wedge P - \aa6 H\wedge M
- \aa{14} P\wedge M \cr
&& \delta(K)= -\frac{\aa2}2 D\wedge P -\aa4 D\wedge M +
\frac{\aa1}2 C \wedge K
+\aa5 C\wedge P \cr
&&\qquad\quad+ \aa8 C\wedge M - \frac{\aa2}2 H\wedge K
-\frac{\aa{11}}3 H\wedge M\cr
&&\qquad\quad +(\aa{10} + \aa4) K\wedge P + \aa{13} K\wedge M -
\frac{\aa{12}}2 P\wedge M \cr
&& \delta(P)= -\frac{ \aa1}2 D \wedge K - \aa3 D \wedge M
-\frac{ \aa1}2 C\wedge P -\frac {\aa7}3 C\wedge M \cr
&&\qquad\quad + \aa5 H\wedge K
+ \frac{\aa2}2 H\wedge P + \aa{10} H\wedge M
- \frac{\aa9 }2 K\wedge M \cr
&&\qquad\quad-(\aa8 + \aa3) K\wedge P + (\aa{13} - \aa6)
P\wedge M \cr
&& \delta(M)=0.
\label{apende}
\eea
 Next, we impose Jacobi identities onto the dual map
$\delta^\ast: \Sch^\ast  \otimes \Sch^\ast \rightarrow \Sch^\ast
$  in order to ensure that $\delta^\ast$ is a  Lie bracket on
$\Sch^\ast$. This condition gives rise to 19 equations
\bea
&&\aa8^2 - \aa3\aa8 - \aa4\aa7 - \aa7\aa{10}/3 = 0 \cr
&&3\aa2\aa{15} - 6\aa4\aa{12} + 4\aa6\aa{11} + 3\aa{10}\aa{12} -
2\aa{11}\aa{13} = 0 \cr
&&3\aa2\aa{10} - 3\aa2\aa4 - 2\aa5\aa{11} = 0 \cr
&&3\aa2\aa6 - 2\aa3\aa{11} - 6\aa4\aa{10} + 3\aa5\aa{12} = 0 \cr
&&\aa3\aa{12} - \aa2\aa{14} - 2 \aa4\aa6 + 2 \aa4\aa{13}
- \aa8\aa{12} + \aa9\aa{11}/3 = 0 \cr
&&3\aa1\aa4 + 2 \aa3\aa5 - 2\aa5\aa8 + \aa2\aa7/3 = 0 \cr
&&\aa1\aa{15} +2\aa4\aa{9}  + 2\aa5\aa{14} -
2\aa8\aa{13} + \aa7\aa{12}/3 = 0\cr
&&2\aa5\aa8 - \aa1\aa{10} - \aa2\aa7 = 0 
\label{apendf}
\eea
\bea
&&\aa{10}^2 - \aa3\aa{11} - \aa4\aa{10} - \aa8\aa{11}/3 = 0  \cr
&&3\aa1\aa{14} - 6\aa3\aa9  - 2\aa6\aa7 - 2\aa7\aa{13}
+ 3\aa8\aa9 = 0  \cr
&&3\aa1\aa3 + 2\aa5\aa7 - 3\aa1\aa8 = 0 \cr
&&3\aa1\aa6 + 6\aa3\aa8 + 2\aa4\aa7 - 3\aa5\aa9 = 0 \cr
&&2\aa3\aa{13} - \aa1\aa{15} + \aa4\aa9 - \aa9\aa{10} +
\aa7\aa{12}/3 = 0 \cr
&&3\aa2\aa3 + 2\aa4\aa5 - 2\aa5\aa{10} + \aa1\aa{11}/3 = 0 \cr
&&\aa2\aa{14} + 2\aa3\aa{12} + 2\aa5\aa{15} +2\aa6\aa{10}
- 2\aa{10}\aa{13} + \aa9\aa{11}/3 = 0 \cr
&&2\aa5\aa{10} - \aa1\aa{11} - \aa2\aa8 = 0 
\label{apendg}
\eea
\bea
&&\aa1\aa2 - \aa5^2 = 0  \cr
&&\aa1\aa{12} + \aa2\aa9 - 2 \aa8\aa{10} + 2\aa7\aa{11}/9 = 0  \cr
&&2\aa3\aa4 + \aa3\aa{10} + \aa4\aa8 - \aa8\aa{10} +
\aa7\aa{11}/9 = 0 . 
\label{apendh}
\eea
Surprisingly, the relations (\ref{apendf}), (\ref{apendg}) and
(\ref{apendh}) can be mapped exactly and in this order  onto the 
 equations (\ref{cb}),  (\ref{cc}) and (\ref{cd}) (which   come
from the modified classical YBE) by means of the following
identification  between the parameters $\aa{l}$ and $\a{i}$,
$\b{i}$, $\c{j}$:
\bea
&&\aa1 = 2\b2  \qquad \aa2 = - 2\a2   \qquad \aa3 = \b1 
\qquad \aa4 = - \a1\qquad \aa5 = - \c3  \cr
&&  \aa6 = - 2 \c1   \qquad \aa7 = - 3 \b5 
\qquad \aa8 = - \a6  \qquad
\aa9 = 2 \b4  \qquad \aa{10} = \b6\label{id}\\
&&   \aa{11} = 3 \a5 
\qquad \aa{12} = - 2 \a4  \qquad
\aa{13} = - (\c1 + \c2)  \qquad \aa{14} = \b3  \qquad
\aa{15} = - \a3.
\nonumber
\eea
 Under this identification, the cocommutator (\ref{apende})
turns into the form (\ref{ci}). Consequently, all   the
Schr\"odinger Lie bialgebras are coboundary ones and
 can be completely described by a classical $r$-marix according
to Theorem 2.1.

\end{document}